\documentclass[12pt,leqno, twoside]{article}

\usepackage{amsfonts}
\usepackage{amsmath}
\usepackage{amssymb}

\topmargin=-0.6cm\usepackage{amssymb}
\usepackage{amsfonts}
\usepackage{amsmath}

\newcommand{\R}{\mathbb{R}}
\newcommand{\C}{\mathbb{C}}
\newcommand{\Oo}{{\cal O}}
\newcommand{\signature}{\operatorname{signature}}
\newcommand{\sgn}{\operatorname{sgn}}
\newcommand{\Aa}{\mathcal{A}}
\newcommand{\K}{\mathbb{K}}
\newcommand{\kropek}{,\ldots ,}
\newcommand{\rank}{\operatorname{rank}}
\newcommand{\ar}{\longrightarrow}
\newcommand{\inv}{^{-1}}

\newtheorem{thm}{Theorem}[section]
\newtheorem{prop}[thm]{Proposition}
\newtheorem{lem}[thm]{Lemma}
\newtheorem{cor}[thm]{Corollary}
\newtheorem{exmp}[thm]{Example}

\numberwithin{equation}{section}

\title{An algebraic formula for the intersection number of a polynomial immersion\thanks{%
							Iwona~Karolkiewicz, Aleksandra~Nowel and Zbigniew~Szafraniec\\
							University of Gda\'{n}sk,
              Institute of Mathematics \\
              80-952 Gda\'{n}sk, Wita Stwosza 57\\
              Poland\\
              Tel.: +48-58-5232059\\
              Fax: +48-58-3414914\\
              Email: ikarolki@manta.univ.gda.pl \\
              Email: Aleksandra.Nowel@math.univ.gda.pl\\
              Email: Zbigniew.Szafraniec@math.univ.gda.pl\\}
}

\author{Iwona~Karolkiewicz \and Aleksandra~Nowel \and Zbigniew~Szafraniec}

\date{July 2007}

\begin{document}

\def\nothanksmarks{\def\thanks##1{\protect\footnotetext[0]{\kern-\bibindent##1}}}
        
\nothanksmarks        

\maketitle

\pagestyle{myheadings}

\markboth{\textsc{\small I.~Karolkiewicz, A.~Nowel, Z.~Szafraniec}}{\textsc{\small Immersions of spheres and algebraically constructible functions}}

\begin{abstract}
There is presented an algorithm for computing the topological degree for
a large class of polynomial mappings.
As an application there is given an effective algebraic formula for
the intersection number of a polynomial immersion
$M \rightarrow \mathbb{R}^{2m}$, where $M$ is an $m$--dimensional
algebraic manifold.
\end{abstract}

{\small

\bigskip


Let $H=(h_1\kropek h_n):\R^n \longrightarrow \R^n$ be a continuous
mapping, and let $p$ be isolated in $H^{-1}(0)$. Then one may
define the local topological degree $\deg_pH$ as the topological
degree of the mapping $$S^{n-1}_r\ni x \longmapsto H(x)/\| H(x)\|
\in S^{n-1},$$ where $S^{n-1}_r$ is the sphere of radius $r\ll 1$
centered at p.

Suppose that $H$ is a polynomial mapping and $u:\R^n
\longrightarrow \R$ is a polynomial. Put $U=\{x\in \R^n:\
u(x)>0\}$. Let $J_{\R}\subset \R[x]=\R[x_1\kropek x_n]$ denote the
ideal generated by $h_1\kropek h_n$, and let $Q=\R[x]/J_{\R}$.
Suppose that $\dim_{\R}Q< \infty$, so that $H^{-1}(0)$ is finite.
According to \cite{Szafraniec}, one may construct quadratic forms
$\Phi_T$ and $\Psi_T$ on Q such that
$$\sum \deg_pH=(\signature \Phi_T +\signature \Psi_T)/2,$$ where $p\in H^{-1}(0)\cap
U$. In \cite{Szafraniec} is also presented a simpler formula for
$$\sum \deg_pH\mod 2.$$

The assumption $\dim_{\R}Q<\infty$ is too restrictive in some
cases. For instance if $H^{-1}(0)$ is infinite then
$\dim_{\R}Q=\infty$.

This paper is devoted to the case where there exists an ideal
$I_{\R}\supset J_{\R}$ such that
$$(J_{\R}:I_{\R})+I_{\R}=\R[x]\quad \mbox{and}\quad
\dim_{\R}\R[x]/(J_{\R}:I_{\R})<\infty.$$(The quotient of ideals
$J_{\R}:I_{\R}$ is defined as the set of all $f\in \R[x]$ with
$fg\in J_{\R}$ for all $g\in I_{\R}$.)

We shall show that in that case $H^{-1}(0)\setminus
V(I_{\R})=V(J_{\R}:I_{\R})$ is finite, and there are effective
algebraic formulae for $\sum \deg_pH\ (\mod 2)$
 and $\sum \deg_pH$ analogous to the ones proved in
 \cite{Szafraniec}.

 Let $M$ be a compact oriented $m$--dimensional manifold ($m>1$),
 and let $g:M \longrightarrow \R^{2m}$ be an immersion. Whitney in
 \cite{whitneySelfInter} introduced the intersection number
 $I(g)$.

Assume that $M\subset \R^{n+m}$ is an algebraic complete
intersection, and a polynomial immersion $g$ has a finite set of
self--intersections. Let $\Delta \subset \R^{n+m}\times  \R^{n+m}$
be the diagonal.

We show that one may construct a polynomial mapping $H: \R^{n+m}
\times  \R^{n+m} \longrightarrow \R^n\times \R^n \times \R^{2m}$
such that $\Delta \subset H^{-1}(0)$, $H^{-1}(0)\setminus \Delta $
is finite, and for even $m$
$$I(g)=\frac{1}{2}\sum_p\deg_pH,$$ where $p\in H^{-1}(0)\setminus
\Delta.$ There is also presented a similar formula for odd $m>1$.

As $\Delta \subset H^{-1}(0) $ is infinite, we cannot apply method
presented in \cite{Szafraniec} so as to compute $\sum_p\deg_pH$.
We show how to apply methods developed in this paper so as to
express the intersection number $I(g)$ either as the signature of
a quadratic form on $\Aa =\R[x]/(J_{\R}:I_{\R})$ (for even $m$),
in terms of $\dim_{\R}\Aa$ and signs of determinants of matrices
of two quadratic forms (for odd $m>1$).

The paper is organized as follows. In Section \ref{idealowo} we
show that the complex algebra $\C[x]/(J_{\C}:I_{\C})$ is
isomorphic to the product of some algebras associated to the ideal
$J_{\C}$ and points in $H_{\C}^{-1}(0)\setminus V(I_{\C})$.  In
Section \ref{part3} we discuss relations between the ideal
$J_{\R}:I_{\R}$ and its complex counterpart $J_{\C}:I_{\C}$. In
Section \ref{part2} we recall some results presented in
\cite{Szafraniec}. In particular we explain how to construct the
bilinear forms which we use in Section \ref{stopien}, where we
prove a formula for $\sum_p\deg_pH$ in terms of their signatures.
In Section \ref{immersje} we apply these results so as to give an
effective method for computing the intersection number of a
polynomial immersion on an algebraic manifold. In the end we
present examples computed by a computer. We implemented our
algorithm with help of {\sc Singular} \cite{GPS06}. We have also
used computer programs written by Adriana Gorzelak and Magdalena
Sarnowska --- students of computer sciences at the University of
Gda\'{n}sk --- which can construct matrices of $\Phi_T$, $\Psi_T$,
and compute the signature of a bilinear form.

\section{Quotients of ideals in $\C[x]$}\label{idealowo}
If $S$ is an ideal in $\C[x]=\C[x_1,\ldots ,x_n]$ then by $V(S)$
we will denote set of all complex zeros of $S$. For each
$p=(p_1,\ldots ,p_n)\in \C^n$, $m_p$ is the maximal ideal
generated by monomials $x_1-p_1, \ldots ,x_n-p_n$. If $J,I\subset
\C[x]$ are ideals, then
$$J:I=\{f\in\C[x]:fg\in J \quad
\textrm{for all}\quad g\in I\}$$
 is called the ideal quotient of
$J$ by $I$. Since $J\subset J:I$, $V(J:I)\subset V(J)$.

Suppose that $J,I\subset \C[x]$ are ideals such that $J\subset I$
and
$$\dim_{\C}\frac{\C[x]}{J:I}<\infty.$$ Then $V(J:I)$ is finite.
Applying \cite[Theorem 7, p. 192]{Cox} we get inclusions:
$$V(J)\setminus V(I)\subset \overline{V(J)\setminus V(I)} \subset
V(J:I)\subset V(J),$$ so that $V(J)\setminus V(I)\subset
V(J:I)\setminus V(I)\subset V(J)\setminus V(I)$. Denote
$V(J)\setminus V(I)=V(J:I)\setminus V(I)= \{ p_1,\ldots , p_s\}$.

\begin{prop}\label{izomor1}
There exists a positive integer $k$ such that
$$\frac{\C[x]}{J} \cong {\frac{\C[x]}{J+I^k}} \times \frac{\C[x]}{J+m_{p_1}^k} \times \ldots \times
\frac{\C[x]}{J+m_{p_s}^k}.$$ In particular $J=(J+I^k)\cap
(J+m_{p_1}^k)\cap \ldots \cap (J+m_{p_s}^k)$. The same holds true
for any integer bigger than $k$.
\end{prop}
{\em Proof.} We have $V(I)\subset V(J)$ and $V(J)\setminus V(I)=
\{ p_1,\ldots , p_s\}$, so $$V(J)=\{p_1,\ldots p_s\}\cup
V(I)=V(m_{p_1}) \cup \ldots \cup V(m_{p_s})\cup V(I)=V(m_{p_1}\cap
\ldots \cap m_{p_s} \cap I).$$ By the Hilbert Nullstellensatz,
there exists $k$ such that \begin{equation}(m_{p_1}\cap \ldots
\cap m_{p_s} \cap I)^k \subset J.\label{inkluzja} \end{equation}
If $i\neq j$ then $V(m_{p_i}^k+m_{p_j}^k)=V(m_{p_i}^k)\cap
V(m_{p_j}^k)=\emptyset$, so,
$$\C[x]=m_{p_i}^k+m_{p_j}^k \subseteq (J+m_{p_i}^k)+(J+m_{p_j}^k)
\subseteq \C[x].$$ and then
$$(J+m_{p_i}^k)+(J+m_{p_j}^k)=\C[x].$$ We
have $V(I^k+m_{p_i}^k)=V(I) \cap \{p_i\}=\emptyset$, so that
$$\C[x]=I^k+m_{p_i}^k\subseteq (J+I^k)+(J+m_{p_i}^k) \subseteq\C[x],$$
and then
$$(J+I^k)+(J+m_{p_i}^k)=\C[x].$$
The Chinese Remainder Theorem implies that:
$$\frac{\C[x]}{(J+I^k)\cap (J+m_{p_1}^k)\cap \ldots \cap (J+m_{p_s}^k)}\cong \frac{\C[x]}{J+I^k} \times \frac{\C[x]}{J+m_{p_1}^k} \times \ldots \times
\frac{\C[x]}{J+m_{p_s}^k}$$ and
$$(J+I^k)\cap (J+m_{p_1}^k)\cap \ldots \cap (J+m_{p_s}^k)=(J+I^k)\cdot (J+m_{p_1}^k)\cdot \ldots \cdot (J+m_{p_s}^k).$$
By (\ref{inkluzja}) we get
$$J\subset (J+I^k)\cap (J+m_{p_1}^k)\cap \ldots \cap (J+m_{p_s}^k)=$$
$$(J+I^k)\cdot (J+m_{p_1}^k)\cdot \ldots \cdot
(J+m_{p_s}^k)\subset  J+(I^{k}\cdot m_{p_1}^k \cdot \ldots \cdot
m_{p_s}^k)=$$
$$J+(I\cdot m_{p_1}\cdot \ldots \cdot
m_{p_s})^k\subset J+(I\cap m_{p_1}\cap \ldots \cap m_{p_s})^k
\subset J+J,$$and then $(J+I^k)\cap (J+m_{p_1}^k)\cap \ldots \cap
(J+m_{p_s}^k)=J$. So we have $$\frac{\C[x]}{J} \cong
{\frac{\C[x]}{J+I^k}} \times \frac{\C[x]}{J+m_{p_1}^k} \times
\ldots \times \frac{\C[x]}{J+m_{p_s}^k}.$$ \hfill $\Box$

\begin{lem}\label{odwracalnosc}
Let  $p\in V(J)$. For any positive integer $k$, $f$ is invertible
in $\C[x]/(J+m_p^k)$ if and only if $f(p)\neq 0$.
\end{lem}
{\em Proof.} ($\Rightarrow$) Since $f$ is invertible in
$\C[x]/(J+m_p^k)$, there exists $g\in \C[x]$ such that $fg\equiv 1
\quad \operatorname{mod}(J+m_p^k)$, so
$$fg=1+w+v,$$ where $w\in J$ and $v\in m_p^k$. Hence
$$f(p)g(p)=1+w(p)+v(p)=1,$$ and so $f(p)\neq 0$.\\($\Leftarrow$) $f$ can be presented as $f=f(p)-g$, where $g\in
m_p$. Then
$$(f(p)-g)(f(p)^{k-1}+f(p)^{k-2}g+ f(p)^{k-3}g^2 +\ldots +g^{k-1})=f(p)^k-g^k.$$
Set $h=(f(p)^{k-1}+f(p)^{k-2}g+\ldots +g^{k-1})$. We have $g^k\in
m_p^k$, so\\ $fh=(f(p)-g)h \equiv f(p)^k \neq 0$ in
$\C[x]/(J+m_p^k)$, and then $f$ is invertible in $C[x]/(J+m_p^k)$.
\hfill $\Box$

Take $k$ as in Proposition \ref{izomor1}. Then
\begin{lem} \label{postac i:j 1}
$J:I=(\bigcap \limits _{i=1}^s(J+m_{p_i}^k))\cap ((J+I^k):I)$
\end{lem}
{\em Proof.} Take $f\in \bigcap \limits _{i=1}^s(J+m_{p_i}^k)\cap
((J+I^k):I)$. For any $g\in I$ and $1\leq i\leq s$ we have
$$fg\in J+I^k \quad \textrm{and} \quad
fg\in J+m_{p_i}^k,$$ so
$fg\in(J+I^k)\cap(J+m_{p_1}^k)\cap\ldots\cap (J+m_{p_s}^k)$. Using
Proposition \ref{izomor1} we get $fg\in J$, so that $f\in J:I$.

Now let us take $f\in J:I$. For any $g\in I$ we have $fg\in
J\subset (J+I^k)$, and then $f\in (J+I^k):I$. For $1\leq i\leq s$
we also have
$$fg \in J+m_{p_i}^k.$$ Because $p_i\notin V(I)$, then there exists
$h \in I$ such that $h(p_i)\neq 0$. From Lemma \ref{odwracalnosc},
$h$ is invertible in $\C[x]/(J+m_{p_i}^k)$. As $fh \equiv 0$ in
$\C[x]/(J+m_{p_i}^k)$ we have
$$f\in J+m_{p_i}^k.$$ So $f\in (\bigcap \limits
_{i=1}^s(J+m_{p_i}^k))\cap ((J+I^k):I)$.
\hfill $\Box$

From now on we shall assume that $(J:I)+I=\C[x]$. We have

\begin{lem}\label{zera}If $(J:I)+I=\C[x]$ and $\dim_{\C} \C[x]/(J:I)<\infty$
then\\ $V(J:I)\cap V(I)=\emptyset$, and then $V(J:I)=V(J)\setminus
V(I)=\{p_1,\ldots,p_s\}$.
\end{lem}
\hfill $\Box$

\begin{lem}
For every positive integer $k$ we have $(J:I)+I^k=\C[x]$.
\end{lem}
{\em Proof.} We have $$V((J:I)+I^k)=V(J:I)\cap V(I^k)=V(J:I)\cap
V(I).$$ Since $(J:I)+I=\C[x]$, we have $V(J:I)\cap
V(I)=\emptyset$. So
$$(J:I)+I^k=\C[x].$$
\hfill $\Box$

\begin{lem}\label{sumaDajeC}
For every positive integer $k$ we have $(J+I^k):I=\C[x]$.
\end{lem}
{\em Proof.} Let us take $f\in (J:I)+I^k$. Then $f=h_1+h_2$, where
$h_1\in J:I$ and $h_2\in I^k$. For any $g\in I$ we get
$$fg=h_1g+h_2g \in J+I^k.$$ So we have  $(J:I)+I^k \subseteq
(J+I^k):I$. Using previous Lemma we get $$\C[x]=(J:I)+I^k
\subseteq (J+I^k):I \subseteq \C[x].$$
\hfill $\Box$

Lemmas \ref{postac i:j 1}, \ref{sumaDajeC} and Proposition
\ref{izomor1} imply

\begin{cor} \label{postac j:i 2}
If $(J:I)+I=\C[x]$ then for all $k$ large enough
$$J:I=\bigcap \limits
_{i=1}^s(J+m_{p_i}^k).$$
\end{cor}
\hfill $\Box$

Ideals $J+m_{p_i}^k$ are pairwise comaximal. As a consequence of
the Chinese Remainder Theorem we get

\begin{prop} \label{izomor2}
If $(J:I)+I=\C[x]$ then for all $k$ large enough
$$\frac{\C[x]}{J:I} \cong \frac{\C[x]}{J+m_{p_1}^k} \times \ldots
\times \frac{\C[x]}{J+m_{p_s}^k}.$$
\end{prop}
\hfill $\Box$


\section{Quotients of ideals in $\R[x]$ and $\C[x]$}\label{part3}
\begin{lem}\label{odRdoC}
Let $f_1,\ldots ,f_r$ be polynomials with real coefficients. Let
$S_{\R}$ (resp. $S_{\C}$) denote the ideal in $\R[x]$ (resp.
$\C[x]$) generated by $f_1,\ldots ,f_r$. Then
\begin{itemize}
\item[(i)] $S_{\R}=S_{\C}\cap \R[x]$, \item
[(ii)]$S_{\R}=\R[x]\Leftrightarrow S_{\C}=\C[x]$.
\end{itemize}
\end{lem}
{\em Proof.} \begin{itemize} \item[(i)] If $f\in S_{\R}$ then
obviously $f\in S_{\C} \cap \R[x]$.

If $h=\sum_{\alpha}a_{\alpha}x^{\alpha}\in \C[x]$ then put
$\overline{h}=\sum_{\alpha}\overline{a_{\alpha}}x^{\alpha}$. In
that case $h\in \R[x]$ if and only if $h=\overline{h}$. Take $f\in
S_{\C} \cap \R[x]$. There exist $h_1, \ldots, h_r \in \C[x]$ such
that $f=\sum \limits _{i=1}^{r} h_i f_i$. As $f, f_1, \ldots ,f_r
\in \R[x]$, we have $\overline{f}=f= \sum \limits _{i=1}^{r}
\overline{h_i} f_i$. Hence $f=\frac{1}{2} \sum \limits _{i=1}^{r}
(h_i+ \overline{h_i} ) f_i$. Of course $h_i+ \overline{h_i} \in
\R[x]$, so $f\in S_{\R}$.

\item[(ii)]If $S_{\R}=\R[x]$ then $1\in S_{\R}\subset S_{\C}$, so
$S_{\C}=\C[x]$.

If $S_{\C}=\C[x]$ then $S_{\R}=S_{\C} \cap \R[x] = \C[x] \cap
\R[x] =\R[x]$.
\end{itemize}
\hfill $\Box$

\begin{lem}
Assume that $J_{\R},I_{\R}$ are ideals in $\R[x]$. Then
$(J_{\R}:I_{\R})_{\C}=J_{\C}:I_{\C}$, and then
$(J_{\C}:I_{\C})\cap \R[x]=J_{\R}:I_{\R}$.
\end{lem}
{\em Proof.} Let $g_1, \ldots, g_s \in \R[x]$ be generators of the
ideal $I_{\R}$, and in consequence of the ideal $I_{\C}$, and let
$f_1, \ldots ,f_l$ generate $J_{\R}$ and also $J_{\C}$. Of course
$$J_{\R}:I_{\R}=\{h\in \R[x]: hg_i \in J_{\R}, \quad \textrm{for each} \quad 1 \leq i \leq s\}$$
$$J_{\C}:I_{\C}=\{h\in \C[x]: hg_i \in J_{\C}, \quad \textrm{for
each} \quad 1 \leq i \leq s\}.$$

Take $h\in (J_{\R}:I_{\R})_{\C}$, then there exist $w_1, \ldots
,w_m \in \C[x]$ and $v_1, \ldots ,v_m \in J_{\R}:I_{\R}$ such that
$$h=\sum_{j=1}^m w_jv_j.$$ For each $1\leq i\leq s$, $v_jg_i\in J_{\R}$, and  $hg_i\in  J_{\C}$, so $h\in J_{\C}:I_{\C}$.

Take $h\in J_{\C}:I_{\C}$, then for each $1\leq i\leq s$, $hg_i
\in J_{\C}$, so there exist $w_1, \ldots ,w_l\in \C[x]$ such that
$hg_i=\sum_{j=1}^l f_jw_j$. Then $\overline{h}g_i=\sum_{j=1}^l
f_j\overline{w}_j$ and
$$(h+\overline{h})g_i=\sum_{j=1}^l(w_j+\overline{w}_j)f_j.$$
Of course $w_j+\overline{w}_j\in \R[x]$. So $(h+\overline{h})g_i
\in J_{\R}$ for each $1\leq i\leq s$, and then $h+\overline{h} \in
J_{\R}:I_{\R} \subset (J_{\R}:I_{\R})_{\C}$. As
$$(h-\overline{h})g_i=\sum_{j=1}^l(w_j-\overline{w}_j)f_j$$ and
$\sqrt{-1}(w_j-\overline{w}_j)\in \R[x]$ so
$\sqrt{-1}(h-\overline{h})g_i\in J_{\R}$ for each $1\leq i\leq s$.
So we get  $\sqrt{-1}(h-\overline{h})\in J_{\R}:I_{\R} $, and then
$h-\overline{h} \in (J_{\R}:I_{\R})_{\C}.$ Because $h+\overline{h}
\in (J_{\R}:I_{\R})_{\C}$ and $h-\overline{h} \in
(J_{\R}:I_{\R})_{\C}$ then $h\in (J_{\R}:I_{\R})_{\C}$.
\hfill $\Box$

If $S_{\K}$ is an ideal in $\K[x]$, where $\K$ is either $\C$ or
$\R$, denote
$$V(S_{\K})=\{p\in \K^n\ :\ f(p)=0\ \mbox{for all}\ f\in S_{\K}\}.$$

Consider ideals $J_{\R}\subset I_{\R} \subset \R[x]$, such that
$$\dim_{\R} \frac{\R[x]}{J_{\R}:I_{\R}} < \infty \quad \textrm{and} \quad
(J_{\R}:I_{\R})+I_{\R}=\R[x].$$ Then
$\C[x]=(J_{\C}:I_{\C})+I_{\C}$ and
$\dim_{\C}\C[x]/(J_{\C}:I_{\C})=\dim_{\R} \R[x]/(J_{\R}:I_{\R}) <
\infty$. By Lemma \ref{zera}, $V(J_{\C}:I_{\C})=V(J_{\C})\setminus
V(I_{\C})$ is finite. Let $V(J_{\C})\setminus V(I_{\C})=\{p_1,
\ldots ,p_s\}$. By Proposition \ref{izomor2} there exists a
positive integer $k$ such that
$$\frac{\C[x]}{(J_{\C}:I_{\C})}  \cong \frac{\C[x]}{J_{\C}+m_{p_1}^k} \times \ldots
\times \frac{\C[x]}{J_{\C}+m_{p_s}^k}.$$

Set $V(J_{\R})\setminus V(I_{\R})=\{p_1,\ldots ,p_m\}$ and
$(V(J_{\C})\setminus V(I_{\C})) \setminus
\R^n=\{q_1,\overline{q_1},\ldots ,q_r,\overline{q_r}\}$. Of course
$s=m+2r$. For $k$ large enough and $p\in V(J_{\R})\setminus
V(I_{\R})$ we define an $\R$--algebra
$$\mathcal{A}_{\R,p}:=\frac{\R[x]}{J_{\R}+m_{\R,p}^k}, \quad \textrm{where}
\quad m_{\R,p}=\{f\in\R[x]:f(p)=0\}.$$ For $p\in
V(J_{\C})\setminus V(I_{\C}) $, we define an $\C$-algebra
$$\mathcal{A}_{\C,p}:=\frac{\C[x]}{J_{\C}+m_{\C,p}^k}, \quad \textrm{where}
\quad m_{\C,p}=\{f\in \C[x]:f(p)=0\}.$$ Of course, $f\in
J_{\C}+m_{\C,p}^k$ if and only if $\overline{f}\in
J_{\C}+m_{\C,\overline{p}}^k$. In particular
\begin{equation} \R[x]\cap (J_{\C}+m_{\C,p}^k)=\R[x]\cap J_{\C}+m_{\C,\overline{p}}^k. \label{wazne1}\end{equation} The mapping $f\longmapsto
\overline{f}$ induces an $\R$--isomorphism of algebras
$\Aa_{\C,p}$ and $\Aa_{\C,\overline{p}}$, so that
\begin{equation} \dim_{\C}\Aa_{\C,p}=2\dim_{\R}\Aa_{\C,p}=2\dim_{\R}\Aa_{\C,\overline{p}}=\dim_{\C}\Aa_{\C,\overline{p}}.\label{wazne2}\end{equation}
 Let us denote
$$\mathcal{B} = \oplus_{i=1}^m \mathcal{A}_{\R,p_i}
\oplus_{j=1}^r \mathcal{A}_{\C,q_j}.$$

\begin{thm}\label{wazny izom} For all $k$ large enough there is a natural isomorphism
$$\frac{\R[x]}{J_{\R}:I_{\R}}\cong \mathcal{B}.$$
\end{thm}
{\em Proof.} We take $k$ as in Proposition \ref{izomor2}. We
define a homomorphism
$$\pi :\R[x]\longrightarrow \mathcal{B},$$
as $\pi (f)=\oplus_{i=1}^m [f]_{p_i} \oplus_{j=1}^r [f]_{q_j}$,
where $[f]_p$ is the residue class of $f$ in the appropriate
algebra. Then
$$\pi^{*} : \frac{\R[x]}{\ker \pi} \longrightarrow \mathcal{B}$$
is a monomorphism. Using Lemma \ref{postac j:i 2} and
(\ref{wazne1}) we get
\begin{align*}
\ker \pi & = \R[x] \cap \bigcap_{i=1}^m (J_{\R}+m_{\R,p_i}^k) \cap
\bigcap_{j=1}^r (J_{\C}+m_{\C,q_j}^k)=\\
& =\R[x] \cap \bigcap_{i=1}^m
(J_{\R}+m_{\R,p_i}^k)_{\C} \cap \bigcap_{j=1}^r
(J_{\C}+m_{\C,q_j}^k)=\\
& =\R[x] \cap \bigcap_{i=1}^m
(J_{\C}+m_{\C,p_i}^k) \cap \bigcap_{j=1}^r
(J_{\C}+m_{\C,q_j}^k)\cap\bigcap_{j=1}^r
(J_{\C}+m_{\C,\overline{q}_j}^k)=\\
& =\R[x]\cap
(J_{\C}:I_{\C})=J_{\R}:I_{\R}.
\end{align*} 
By (\ref{wazne2}) we have
\begin{align*}
\dim_{\R} \mathcal{B} & =\sum_{i=1}^m \dim_{\R} \mathcal{A}_{\R,p_i}
+\sum_{j=1}^r \dim_{\R} \mathcal{A}_{\C,q_j}=\\
& =\sum_{i=1}^m
\dim_{\C} \mathcal{A}_{\C,p_i} +2\cdot \sum_{j=1}^r \dim_{\C}
\mathcal{A}_{\C,q_j}=\\
& =\sum_{i=1}^m \dim_{\C}
\mathcal{A}_{\C,p_i} +\sum_{j=1}^r \dim_{\C}
\mathcal{A}_{\C,q_j}+\sum_{j=1}^r \dim_{\C}
\mathcal{A}_{\C,\overline{q}_j}\ .
\end{align*}

By Proposition \ref{izomor2}, it equals
$$\dim_{\C}\frac{\C[x]}{J_{\C}:I_{\C}}=\dim_{\R}
\frac{\R[x]}{J_{\R}:I_{\R}}.$$ So $\dim_{\R} \mathcal{B} =
\dim_{\R} \R[x]/(J_{\R}:I_{\R})$ and then
$\pi^*:\R[x]/(J_{\R}:I_{\R})\rightarrow \mathcal{B}$ is an
isomorphism.
\hfill $\Box$


\section{Bilinear forms}\label{part2}
Let $\K$ denote either $\R$ or $\C$. For $p\in \K^n$, let $\Oo
_{\K,p}$ denote the ring of germs at $p$ of analytic functions
$\K^n\longrightarrow \K$. There is a natural homomorphism
$$\eta :\K[x]=\K[x_1, \ldots ,x_n]\longrightarrow \Oo_{\K,p}.$$ Let
$m_{\K,p}=\{f\in \K[x]:\  f(p)=0\}$ be the maximal ideal in
$\K[x]$ associated with $p$.

Let $S_{\R}$ be an ideal in $\R[x]$, let $S_{\C}$ denote the ideal
in $\C[x]$ generated by $S_{\R}$, and let
$$V(S_{\K})=\{p\in \K^n:\ f(p)=0 \ \mbox{for all } f\in S_{\K}\}.$$

Take $h_1, \ldots ,h_n \in S_{\R}$, the ideal $J_{\K}$ generated
by $h_1\kropek h_n$, and a polynomial mapping $H_{\K}=(h_1,\ldots
,h_n):\K^n \longrightarrow \K^n$. Then $J_{\K}\subset S_{\K}$ and
$V(S_{\K})\subset V(J_{\K})=H_{\K}^{-1}(0)$. In particular, points
isolated in $H_{\C}^{-1}(0)$ are isolated in $V(S_{\K})$.

Take $p\in V(S_{\K})$. Let $S_{\K,p}$ (resp. $J_{\K,p}$) denote
the ideal in $\Oo_{\K,p}$ generated by $S_{\K}$ (resp. $J_{\K}$),
and let $\Aa'_{\K,p}=\Oo_{\K,p} / S_{\K,p}$. Clearly $\Aa'_{\K,p}$
is an $\K$--algebra and $\eta (S_{\K})\subset S_{\K,p}$.

Assume that $S_{\K,p}=J_{\K,p}$ so that $\Aa'_{\K,p}=\Oo_{\K,p} /
J_{\K,p}$. Then $\dim_{\K} \Aa'_{\K,p} < \infty$ if and only if
$p$ is isolated in $H_{\C}^{-1}(0)$. If that is the case then $p$
is isolated in $V(S_{\R})$, and $\eta (m^k_{\K,p}) \subset
S_{\K,p}$ for all $k$ large enough.
\begin{lem}\label{pi indukuje}
If $p\in V(S_{\K})$ is isolated in $H_{\C}^{-1}(0)$ and
$S_{\K,p}=J_{\K,p}$, then $\eta$ induces an isomorphism of
$\K$--algebras $$\eta: \K[x] / (S_{\K}+m^k_{\K,p}) \longrightarrow
\Aa'_{\K,p}$$ for all $k$ large enough.
\end{lem}
\hfill $\Box$

For $p$ and $k$ as above, put $$\Aa_{\K,p}=\K[x]/
(S_{\K}+m_{\K,p}^k),$$ so that $\Aa_{\K,p}$ is isomorphic to
$\Aa'_{\K,p}$. In particular, $\Aa_{\K,p}$ does not depend on $k$,
if $k$ is large enough, and $\dim_{\K} \Aa_{\K,p} < \infty$.

Applying the formula for the local topological degree by Eisenbud
and Levine \cite{E&L} and Khimshiashvili \cite{khims1, khims2} ,
and the theory of Frobenius algebras (see \cite{BCRSz, Kunz, s_s})
one may prove some properties, presented in \cite{Szafraniec}, of
bilinear forms on the algebra $\Oo_{\K,p} / J_{\K,p}=\Aa'_{\K,p}$.
As $\Aa'_{\K,p}$ is isomorphic to $\Aa_{\K,p}$, they hold true for
the algebra $\Aa_{\K,p}$.

Let points $p_1, \ldots ,p_w \in V(S_{\C})$ be as in Lemma \ref{pi
indukuje}. One may assume that $p_1, \ldots ,p_m \in \R^n$ and
$p_{m+1}, \ldots , p_w \in \C^n \setminus \R^n$. Denote
$${\mathcal B}=\oplus_{i=1}^m \Aa_{\R,p_i} \oplus _{j=m+1}^w
\Aa_{\C,p_j}.$$ Obviously, ${\mathcal B}$ is a finite dimensional
$\R$--algebra.

Let $u\in\R[x]$, and let $\varphi:{\mathcal B} \longrightarrow \R$
be a linear functional. Then there are bilinear symmetric forms
 $\Phi,\,\Psi:{\mathcal B} \times{\mathcal B}\longrightarrow \R$ given by
$\Phi(a,b)=\varphi(ab)$ and $\Psi(a,b)=\varphi(uab)$.

One may define $\operatorname{signature}\Phi$, and similarly
$\operatorname{signature}\Psi$, as the dimension of a maximal
subspace of ${\mathcal B}$ on which $\Phi$ is positive definite
minus the dimension of that one on which $\Phi$ is negative
definite.

Let $\det [\Phi]$ (resp. $\det [\Psi]$) denote the determinant of
the matrix of $\Phi$ (resp. $\Psi$), with respect to some basis of ${\mathcal
B}$. The sign of the determinant does not depend on the choice of
a basis.

\begin{thm}{\rm \cite[Theorem 2.3, p. 306]{Szafraniec}} \label{zyrandol} Let points $p_1, \ldots ,p_w \in V(S_{\C})\subset H_{\C}^{-1}(0)$ be as in Lemma \ref{pi indukuje}. Suppose that $\det[\Psi]\neq0$. Then
$\det[\Phi]\neq 0$, $u(p_i)\neq 0$ for $1\leq i\leq m$, and
\[\sum\deg_p H_{\R}\equiv\dim_{\R}{\mathcal B}+1+(\sgn \det[\Phi]+\sgn \det[\Psi])/2
\pmod{2},\] where $p\in\{p_1,\ldots,p_m\}\cap \{u>0\}$ and $\deg_p
H_{\R}$ denotes the local topological degree of $H_{\R}$ at
$p$.\end{thm}
\hfill $\Box$

For $x=(x_1,\ldots ,x_n),\ x'=(x_1',\ldots ,x_n')$, and $1\leq
i,j\leq n$ define
\[ T_{ij}(x,x')=\frac{h_i(x_1',\ldots ,x_{j-1}',x_j, \ldots , x_n)-
   h_i(x_1',\ldots ,x_j',x_{j+1},\ldots ,x_n)}{x_j-x_j'}.\]
It is easy to see that each $T_{ij}$ extends to a polynomial, thus
we may assume that
\[T_{ij}\in\R[x,x']=\R[x_1,\ldots ,x_n,x_1',\ldots ,x_n'].\]
There is the natural projection $\R[x,x']\longrightarrow {\mathcal
B}\otimes{\mathcal B}$ given by
\[x_1^{\alpha_1}\cdots x_n^{\alpha_n}(x_1')^{\beta_1}\cdots (x_n')^{\beta_n}
\mapsto x_1^{\alpha_1}\cdots x_n^{\alpha_n}\otimes
         (x_1')^{\beta_1}\cdots (x_n')^{\beta_n}.\]
Let $T$ denote the image of $\det [T_{ij}(x,x')]$ in ${\mathcal
B}\otimes{\mathcal B}$.

Put $d=\dim_{\R}{\mathcal B}$. Assume that $e_1,\ldots ,e_d$ form
a basis in ${\mathcal B}$. So $\dim_{\R} {\mathcal B}\otimes
{\mathcal B}=d^2$ and $e_i\otimes e_j$, for $1\leq i,j\leq d$,
form a basis in ${\mathcal B}\otimes{\mathcal B}$. Hence there are
$t_{ij}\in\R$ such that
\[ T=\sum_{i,j=1}^d t_{ij}e_i\otimes e_j=\sum_{i=1}^d e_i\otimes\hat{e}_i,\]
where $\hat{e}_i=\sum_{j=1}^d t_{ij}e_j.$ Elements
$\hat{e}_1,\ldots ,\hat{e}_d$ form a basis in ${\mathcal B}.$ So
there are $A_1,\ldots , A_d\in\R$ such that
\[1=A_1\hat e_1+\ldots  +A_d\hat e_d\mbox{ in }{\mathcal B}.\\[1em]\]
{\em Definition.\/} For $f=a_1e_1+\ldots + a_de_d\in{\mathcal B}$
define $\varphi_T(f)=a_1A_1+\ldots +a_d A_d.$ Hence
$\varphi_T:{\mathcal B}\longrightarrow \R$
is a linear functional.\\[1em]

Let $\Phi_T$ be the bilinear form on ${\mathcal B}$ given by
$\Phi_T(a,b)=\varphi_T(ab)$.

\begin{thm}{\rm \cite[Theorem 1.5, p. 304]{Szafraniec}}\label{wklad}
 The form $\Phi_T$ is non-degenerate and
$$\sum_{i=1}^m \deg_{p_{i}} H_{\R}=\operatorname{signature}
\Phi_T.$$
\end{thm}
\hfill $\Box$


\section{Topological degree}\label{stopien}

Let $h_1,\ldots,h_n\in\R[x_1,\ldots,x_n]$, and let
$H_{\K}=(h_1,\ldots,h_n):\K^n\longrightarrow \K^n$. Denote by
$J_{\K}$ the ideal in $\K[x]$ generated by $h_1,\ldots,h_n$, so
that $H_{\C}^{-1}(0)=V(J_{\C})$.

Assume that there is an ideal $I_{\R}$ such that $J_{\C}\subset
I_{\C}$, $\dim_{\R}\R[x]/(J_{\R}:I_{\R})<\infty$, and
$(J_{\R}:I_{\R})+I_{\R}=\R[x]$.

Put $S_{\R}=J_{\R}:I_{\R}$. Hence $V(S_{\C})$ is finite, and by
Lemma \ref{zera}
$$H_{\C}^{-1}(0)\setminus V(I_{\C})=V(J_{\C})\setminus V(I_{\C})= V(S_{\C})\ .$$

Hence each $p\in V(J_{\C})\setminus V(I_{\C})$ is isolated in
$H_{\C}^{-1}(0)$. By Corollary \ref{postac j:i 2}, if $k$ is large
enough then the ideal $S_{\K,p}\subset {\cal O}_{\K,p}$ generated
by $S_{\R}$ equals the ideal generated by $J_{\K}+m_{\K,p}^k$.

Let $J_{\K,p}$ denote the ideal in ${\cal O}_{\K,p}$ generated by
$J_{\R}$. Since $J_{\K,p}$ has an algebraically isolated zero at
$p$, the local Nullstellensatz implies that $m_{\K ,p}^k\subset
J_{\K ,p}$, so that $S_{\K,p}=J_{\K ,p}+m_{\K ,p}^k=J_{\K ,p}$.
Hence each point $p\in V(S_{\C})=H_{\C}^{-1}(0)\setminus
V(I_{\C})$ satisfies the assumptions of Lemma \ref{pi indukuje}.

Put $H_{\R }^{-1}(0)\setminus V(I_{\R })=\{p_1,\ldots,p_m\}$ and
$$(H_{\C }^{-1}(0)\setminus V(I_{\C }))\setminus\R^n=
\{q_1,\overline{q_1},\ldots,q_r,\overline{q_r}\}\ .$$

By Theorem \ref{wazny izom}, $\Aa:=\R[x]/S=\R[x]/(J_{\R}:I_{\R})$
and
$${\mathcal B}=\oplus_{i=1}^m \Aa_{\R,p_i} \oplus _{j=1}^r
\Aa_{\C,q_j}$$ are isomorphic.

As a consequence of Theorem \ref{zyrandol} we get

\begin{thm}\label{teczka} Assume that
$\dim_{\R}\Aa <\infty$, $u\in\R[x]$ and $\varphi
:\Aa\rightarrow\R$ is a linear functional. Let bilinear symmetric
forms
 $\Phi,\,\Psi:{\mathcal A} \times{\mathcal A}\longrightarrow \R$ be given by
$\Phi(a,b)=\varphi(ab)$ and $\Psi(a,b)=\varphi(uab)$.

Suppose that $\det[\Psi]\neq0$. Then $\det[\Phi]\neq 0$, $u(p)\neq
0$ for each $p\in H_{\R}^{-1}(0)\setminus V(I_{\R})$, and
\[\sum\deg_p H_{\R}\equiv\dim_R \Aa+1+(\sgn \det[\Phi]+\sgn \det[\Psi])/2
\pmod{2},\] where $p\in H_{\R}^{-1}(0)\cap \{u>0\}\setminus
V(I_{\R})$. \end{thm}
\hfill $\Box$

The same way as in Section  \ref{part2} one may construct the
bilinear symmetric form $\Phi_T$ on $\Aa=\R[x]/(J_{\R}:I_{\R})$.
As a consequence of Theorem \ref{wklad} we get

\begin{thm}\label{drzewo}
The form $\Phi_T$ is non-degenerate and
$$\signature \Phi_T=\sum \deg_p H_{\R}\ ,$$
where $p\in H_{\R}^{-1}(0)\setminus V(I_{\R})$.
\end{thm}
\hfill $\Box$

Let $\Psi_T$ be the bilinear form on $\Aa$ given by
$\Psi_T(a,b)=\varphi_T(uab)$. Using the same arguments as in
\cite[Theorem 1.5, p. 304]{Szafraniec} one may prove
\begin{thm} $\Psi_T$ is non-degenerate if and only if $u(p)\neq
0$ for each $p\in H_{\C}^{-1}(0)\setminus V(I_{\C})$. If that is
the case then $$\sum \deg_pH_{\R}=\frac{1}{2}(\signature \Phi_T
+\signature \Psi_T),$$ where $p\in H_{\R}^{-1}(0)\cap
\{u>0\}\setminus V(I_{\R})$.
\end{thm}

\section{Immersions}\label{immersje}

Let $M$ be an $m$--dimensional manifold. A $C^1$ map $g:M
\longrightarrow \R^k$ is called an immersion, if for each $p\in M$
the rank of $Dg(p)$ equals $m$.

Assume that $m$ is even, $k=2m$ and $M$ is compact and oriented.
Define $$G:M\times M \longrightarrow \R^{2m} \quad \textrm{as}
\quad G(x,y)=g(x)-g(y).$$ Set $\Delta =\{(p,p):p\in M\} \subset
M\times M$. Since $g$ is an immersion, $\Delta$ is isolated in
$G^{-1}(0)$, and so $G^{-1}(0) \setminus \Delta$ is a compact
subset of $M\times M \setminus \Delta$. There exists $(N, \partial
N)$ --- a compact $2m$-dimensional oriented manifold with boundary,
such that $$N\subset M\times M \setminus \Delta \quad \textrm{and}
\quad G^{-1}(0) \setminus \Delta \subset N\setminus
\partial N.$$

Denote by $d(G)$ the topological degree of the mapping $$\partial
N\ni (x,y)\mapsto \frac{G(x,y)}{\| G(x,y)\| } \in S^{2m-1}.$$ Of
course, $d(G)$ does not depend on the choice of $N$. In
particular, if $G^{-1}(0)\setminus \Delta$ is finite then
$$d(G)=\sum \deg_zG, \quad \textrm{where} \quad z\in G^{-1}(0)\setminus
\Delta,$$and $\deg_zG$ denotes the local topological degree of $G$
at $z$.

Whitney has introduced in \cite{whitneySelfInter} an intersection
number $I(g)$ for the immersion $g$. According to the Lashof and
Smale \cite[Theorem 3.1]{LashofSmale},
$$d(G)=2I(g).$$

Now assume that $f=(f_1,\ldots \,f_n):\R^{n+m} \longrightarrow
\R^n$ is a $C^1$ mapping, such that $M=f^{-1}(0)$ and $M$ is a
complete intersection, i.e. for each $p\in M$ the rank of $Df(p)$
equals $n$.

We shall say that vectors $v_1, \ldots ,v_m\in T_pM$ are well
oriented if vectors \\$\nabla f_1(p),\ldots ,\nabla
f_n(p),v_1,\ldots ,v_m$ are well oriented in $\R^{n+m}$. Put
$$F(x,y)=(f(x),f(y)):\R^{n+m}\times \R^{n+m} \longrightarrow
\R^n\times \R^n.$$ Then $F^{-1}(0)=M\times M$ is a complete
intersection. The orientation of $F^{-1}(0)$ defined the way as
above is the same as the orientation of the product $M\times M$.

Let $$\overline{g}=(g_1,\ldots ,g_{2m}):\R^{n+m} \longrightarrow
\R^{2m}$$ be a $C^1$ mapping. Put $g=\overline{g}|_M$. It easy to
verify that $$\rank \left [ \begin{array}{c} Dg(p)\\
\end{array} \right ] = \rank \left [
\begin{array}{c}
  D\overline{g}(p)\\
  Df(p)\\
\end{array} \right ]-n$$
at each point $p\in M$, so we have
\begin{lem}\label{noga}
$g=\overline{g}|_M:M\longrightarrow \R^{2m}$ is an immersion if
and only if at each point $p\in M$ $$\rank \left [
\begin{array}{c}
  D\overline{g}(p)\\
  Df(p)\\
\end{array} \right ]=n+m.$$
\end{lem}
\hfill $\Box$
\begin{exmp}\label{przyklad1}
Let $f=x^2+y^2+z^2-r^2$, and
$$\overline{g}=(x,y,xz,yz):\R^3 \longrightarrow \R^4,$$ i.e. $m=2$,
$n=1$. Then $M=f^{-1}(0)$ is the $2-$dimensional sphere $S^2(r)$
of radius $r$. As $$\rank \left [
\begin{array}{c}
  D\overline{g}\\
  Df\\
\end{array} \right ]=\rank \left [
\begin{array}{ccc}
  1 & 0 & 0\\
  0 & 1 & 0\\
  z & 0 & x\\
  0 & z & y\\
  2x & 2y & 2z\\
\end{array} \right ]$$
has a non-zero $(3\times 3)$--minor at each point $p\in
\R^3\setminus \{0\}$, then $g=\overline{g}|_{S^2(r)}$ is an
immersion for every $r>0$.
\end{exmp}

Let $$\overline{G}(x,y)=\overline{g}(x)-\overline{g}(y).$$ Then
$G(x,y)=g(x)-g(y)=\overline{G}|_{M\times
M}=\overline{G}|_{F^{-1}(0)}$. Put
$$H(x,y)=(F(x,y),\overline{G}(x,y))=$$$$=(f_1(x),\ldots ,f_n(x),f_1(y),\ldots
,f_n(y),g_1(x)-g_1(y),\ldots ,g_{2m}(x)-g_{2m}(y)).$$ Then
$H:\R^{n+m}\times \R^{n+m}\longrightarrow \R^{2n+2m}$, and
$(p,q)\in H^{-1}(0)$ if and only if $(p,q)\in
M\times M$ and $G(p,q)=0$. By \cite[Lemma
3.2]{Szafraniec}, $z=(p,q)\in M\times M$ is isolated in
$G^{-1}(0)$ if and only if $z$ is isolated in $H^{-1}(0)$, and if
that is the case then
$$\deg_zG=\deg_zH.$$
Let $\Delta =\{(p,p)|p\in \R^{n+m}\}$ be the diagonal in
$\R^{n+m}\times \R^{n+m}$. Then $$z=(p,q)\in H^{-1}(0)\setminus
\Delta \quad \textrm{if and only if} \quad z=(p,q)\in
G^{-1}(0)\setminus \Delta.$$ So we get
\begin{prop}\label{stopien+imm}
Suppose that $m$ is even, and
\begin{itemize}
\item[(a)]$M=f^{-1}(0)$ is an oriented compact $m$--dimensional
complete intersection,
\item[(b)]$g=\overline{g}|_M:M\longrightarrow \R^{2m}$ is an
immersion, \item[(c)]$H^{-1}(0)\setminus \Delta$ is finite.
\end{itemize}
Then $2I(g)=d(G)=\Sigma \deg_zH$, where  $z\in H^{-1}(0)\setminus
\Delta.$
\end{prop}
\hfill $\Box$

A homotopy $h_t\colon M\ar \R ^{2m}$ is called \emph{a regular
homotopy}, if at each stage it is an immersion and the induced
homotopy of the tangent bundle is continuous.

As in \cite{whitneySelfInter} we say that an immersion $g\colon
M\ar \R^{2m}$ has \emph{a regular self--intersection} at the point
$g(p)=g(q)$ if
$$D g(p)T_pM+D g(q)T_qM=\R ^{2m}.$$
That is so if and only if $\det [DH(p,q)]\neq 0$.

If $g$ has only regular self--intersections, and there are no
triple points $g(p)=g(q)=g(w)$, then we say that $g$ is completely
regular. If $m$ is odd, then one can define \emph{the intersection
number} of a completely regular immersion as the number of its
self--intersections modulo $2$.

\medskip

By \cite[Theorem 2]{whitneySelfInter}, if $M$ is compact then the
intersection number is invariant under regular homotopies. As in
\cite[\S 8]{whitneyDiffMan}, for any immersion $g\colon M\ar
\R^{2m}$ there exists a completely regular immersion
$\tilde{g}\colon M\ar \R ^{2m}$ which is arbitrarily close to $g$
in $C^1$--topology and there is a regular homotopy between $g$ and
$\tilde{g}$. Thus if $m>1$ is odd then one can define the
intersection number $I(g)$ as the number of self--intersections of
$\tilde{g}$ modulo $2$.

Suppose that $g$ has a finite number of self--intersections, i.e.
$H\inv(0)\setminus \Delta$ is finite. The immersion has a
self--intersection $g(p)=g(q)$ if and only if $(p,q)$ and $(q,p)$
belong to $H\inv (0)\setminus \Delta$. There exists a linear
function
$$u(x,y)=a_1(x_1-y_1) +\ldots +a_{n+m}(x_{n+m}-y_{n+m})$$
which does not vanish at any point in $H\inv (0)\setminus \Delta$,
so that $u(p,q)>0$ if and only if $u(q,p)<0$. Then each
self--intersection is represented by a single point in $H\inv
(0)\cap \{ u>0\} \setminus \Delta$.

\begin{prop}\label{odd}
Suppose that $m>1$ is odd, and
\begin{itemize}
\item[(a)] $M=f\inv (0)$ is an oriented, compact $m$--dimensional
complete intersection, \item[(b)] $g=\overline{g}|M\colon M\ar \R
^{2m}$ is an immersion, \item[(c)] $H\inv (0)\setminus \Delta$ is
finite, \item[(d)] $u(x,y)=a_1(x_1-y_1) +\ldots
+a_{n+m}(x_{n+m}-y_{n+m})$ is such that $u(p,q)\neq 0$ for each
$(p,q)\in H\inv (0)\setminus \Delta$.
\end{itemize}
Then
$$I(g)\equiv \sum \deg _zH \mod 2,$$
where $z\in H\inv (0)\cap\{ u>0\} \setminus \Delta$.
\end{prop}

{\it Proof.} Put
$$\tilde{H}(x,y)=$$$$=(f_1(x),\ldots,f_n(x),f_1(y),\ldots,f_n(y),\tilde{g}_1(x)-\tilde{g}_1(y),\ldots,\tilde{g}_{2m}(x)-\tilde{g}_{2m}(y)).
$$ Let $B_z\subset \R^{n+m}\times \R^{n+m}$ denote a small ball
centered at $z\in H\inv (0)\cap\{ u>0\} \setminus \Delta$. If
$\tilde{g}$ is close enough to $g$, then self--intersections of
$\tilde{g}$ are represented by points in the set $\tilde{H}\inv
(0)\cap\{ u>0\} \setminus \Delta$, which is a subset of the union
of all $B_z$. As self--intersections of $\tilde{g}$ are regular,
$$I(\tilde{g})\equiv \sum _w \sgn \det [D\tilde{H}(w)] \mod 2,$$
where $w\in \tilde{H}\inv (0)\cap\{ u>0\} \setminus \Delta$. Each
point $w$ belongs to some $B_z$, hence
$$I(g)=I(\tilde{g})\equiv \sum _z\sum _w \sgn \det [D\tilde{H}(w)] \mod 2,$$
where $z\in H\inv (0)\cap\{ u>0\} \setminus \Delta$ and $w\in
\tilde{H}\inv (0)\cap B_z$. Since $\tilde{H}$ is close to $H$ in a
neighbourhood of $M\times M$,
$$I(g)\equiv \sum _z \deg _z H \mod 2.$$
\hfill $\Box$

If $f_1\kropek f_n,g_1\kropek g_{2m}$ are polynomials then
$H=(h_1\kropek h_{2n+2m})$ is a polynomial mapping. Let $J_{\R}$
denote the ideal in $\R[x,y]=\R[x_1\kropek x_{n+m}, y_1\kropek
y_{n+m}]$ generated by $h_1\kropek h_{2n+2m}$, and $I_{\R}$ the
one generated by $$f_1(x)\kropek f_n(x),f_1(y)\kropek
f_n(y),x_1-y_1\kropek x_{n+m}-y_{n+m}.$$ It is easy to verify that
$J_{\R}\subset I_{\R}$. Then $$V(J_{\R})=H^{-1}(0),\
V(I_{\R})=M\times M\cap \Delta.$$ Let $Q=\R[x,y]/J_{\R}$. If
$m\geq 1$ and $M\neq \emptyset$, then $M\times M\cap \Delta
\subset H^{-1}(0)$ is infinite, so $\dim_{\R}Q=\infty$ and we
cannot apply methods presented in \cite{Szafraniec} so as to
compute $\sum_p\deg_pH$.

Let $\Aa=\R[x,y]/(J_{\R}:I_{\R}).$ Suppose that $\dim_{\R}\Aa
<\infty$ and $(J_{\R}:I_{\R})+I_{\R}=\R[x,y]$. Then
$H^{-1}(0)\setminus \Delta$ is finite, so that the immersion $g$
has a finite set of self--intersections.

Let $\Phi_T :\Aa\times \Aa\longrightarrow \R$ be the bilinear form
constructed as in Section \ref{part2}. As a consequence of Theorem
\ref{drzewo} and Proposition \ref{stopien+imm} we get
\begin{thm}\label{nowe1}
If $m$ is even, then $$I(g)=\frac{1}{2} \signature \Phi_T.$$
\end{thm}
\hfill $\Box$

Take any polynomial $u(x,y)=a_1(x_1-y_1)+\ldots
+a_{n+m}(x_{n+m}-y_{n-m})$ and any linear from $\varphi:\Aa
\longrightarrow \R$. Let $\Phi ,\Psi:\Aa\times \Aa \longrightarrow
\R$ be bilinear forms given by $\Phi(a,b)=\varphi(ab)$,
$\Psi(a,b)=\varphi(uab)$. As a consequence of Theorem \ref{teczka}
and Proposition \ref{odd} we get
\begin{thm}\label{nowe2}
If $m>1$ is odd and $\det[\Psi]\neq 0$, then $$I(g)\equiv
\dim_{\R}\Aa +1+(\sgn \det[\Phi]+\sgn \det[\Psi])/2\mod2.$$
\end{thm}
\hfill $\Box$


\begin{exmp}
Let us consider the mapping $$g=(g_1,g_2,g_3,g_4)=(x_1, x_2,
x_1x_3,x_2x_3):\R^3 \longrightarrow \R^4.$$ As in Example
\ref{przyklad1}, $g|_{S^2(1)}$ is an immersion. Put $f
=x_1^2+x_2^2+x_3^2-1$. With the  immersion $g|_{S^2(1)}$ we may
associate the polynomial mapping $$H=(h_1\kropek h_6):\R^3\times
\R^3 \longrightarrow \R^6$$ given by $h_1=f(x_1,x_2,x_3)$,
$h_2=f(y_1,y_2,y_3)$ and
$h_i=g_{i-2}(x_1,x_2,x_3)-g_{i-2}(y_1,y_2,y_3)$, for $i=3,4,5,6$.

Let $I_{\R}$ be the ideal in $\R[x_1,x_2,x_3,y_1,y_2,y_3]$
generated by $f(x),f(y),x_1-y_1,x_2-y_2,x_3-y_3$ and let $J_{\R}$
be generated by $h_1\kropek h_6$. Using {\sc Singular} --- a
computer algebra system for polynomial computations --- one may
check that $J_{\R}:I_{\R}$ is generated by
$x_1,x_2,y_1,y_2,x_3+y_3,y_3^2-1$, and then monomials $e_1=1$ and
$e_2=y_3$ form a basis of $\Aa =\R[x_1\kropek
y_3]/(J_{\R}:I_{\R})$. Then
 $\dim_{\R}\Aa=2<\infty$ and $(J_{\R}:I_{\R})+I_{\R}=\R[x_1,x_2,x_3,y_1,y_2,y_3]$.
One may check that
$$T=-8y_3y_3'-8 \quad \mbox{in}\  \Aa \otimes \Aa.$$ So $\hat{e}_1=-8$ and
$\hat{e}_2=-8y_3'$. Clearly $1=(-\frac{1}{8})\hat{e}_1+0\hat{e}_2$
in $\Aa$. Then $A_1=-\frac{1}{8}$, $A_2=0$, so for any
$a=a_1+a_2y_3=a_1e_1+a_2e_2$ in $\Aa$, $\varphi_T$ is given by
$\varphi_T(a)=-\frac{a_1}{8}$. Then the matrix of $\Phi_T$ is
given by $$\left[ \begin{array}{cc} -\frac{1}{8} & 0
\\0 & -\frac{1}{8} \end{array} \right],$$ so $\signature \Phi_T=-2$, and as a consequence of Theorem \ref{nowe1} we get $I(g|_{S^2(1)})=-1$.
\end{exmp}

\begin{exmp}
Let $g=(g_1\kropek g_6)=(x_1,x_2,x_1x_3,x_2x_3,x_4,x_3x_4):\R^4
\longrightarrow \R^6$. Using Lemma \ref{noga} it is easy to check
that $g|_{S^3(1)}$ is an immersion. Put
$f=x_1^2+x_2^2+x_3^2+x_4^2-1$ and
$$H(x_1\kropek y_4)=(h_1\kropek h_8)=$$$$=(f(x_1\kropek x_4),f(y_1\kropek
y_4),g(x_1\kropek x_4)-g(y_1\kropek y_4)).$$ Let $I_{\R}\subset
\R[x_1\kropek y_4]$ be the ideal generated by
$f(x),f(y),x_1-y_1,x_2-y_2,x_3-y_3,x_4-y_4$, and $J_{\R}$
generated by $h_1\kropek h_8 $.  One may check, using {\sc
Singular}, that $x_1,x_2,x_4,y_1,y_2,y_4,x_3+y_3,y_3^2-1$ generate
$J_{\R}:I_{\R}$, and that monomials $e_1=1,\ e_2=y_3$ form the
basis of $\Aa=\R[x_1\kropek y_4]/(J_{\R}:I_{\R})$. Then
 $\dim_{\R}\Aa =2<\infty$ and $(J_{\R}:I_{\R})+I_{\R}=\R[x_1\kropek y_4]$. Put
$\varphi :\Aa \longrightarrow \R$ as
$$\varphi (a)=\varphi (a_1e_1+a_2e_2)=a_2.$$Take
$u=x_3-y_3$. Matrices of $\Phi$ and $\Psi$ are given by

$$\left[
\begin{array}{cc} 0 & +1
\\+1 & 0 \end{array} \right],\quad \left[ \begin{array}{cc} -2& 0
\\0 & -2 \end{array} \right].$$
Of course $\det[\Psi]\neq 0$, and as consequence of  Theorem
\ref{nowe2} we get $I(g|_{S^3(1)})\equiv 2+1+\frac{1}{2}(-1+1)
\equiv 1\ \mod 2 $.
\end{exmp}

Using similar methods we have computed some more difficult
examples:
\begin{exmp}
$h(x_1,x_2,x_3)=$$$=(2x_1x_2+x_2,2x_1x_3+4x_3,4x_3^2+5x_2,5x_2^2+4x_3)$$
is an immersion on the $2-$dimensional sphere of radius $r=10$. In
that case $\dim_{\R}\Aa= 16$, and $I(h|_{S^2(10)})=0$.\end{exmp}
\begin{exmp}
$h(x_1,x_2,x_3)=$$$=(5x_2x_3+x_3^2+3x_1,4x_1^2+3x_3^2+x_3,2x_2^2+3x_2x_3+2x_1,x_2x_3+4x_3^2+3x_2)$$
is an immersion on the $2-$dimensional spheres of radius $r=1,10$.
In both cases $\dim_{\R}\Aa= 6$, and $I(h|_{S^2(1)})=0$,
$I(h|_{S^2(10)})=1$.\end{exmp}
\begin{exmp}$h(x_1,x_2,x_3)=$$$=(3x_1x_2+2x_2^2+2x_1x_3+3x_1+5x_3,2x_1x_2+5x_2^2+3x_2x_3+x_1+2x_2,$$$$4x_1^2+4x_1x_3+5x_2x_3+3x_1+3x_3,
4x_2^2+3x_1x_3+4x_2x_3+4x_1+4x_3)$$
 is an immersion
on the $2-$dimensional spheres of radius $r=1$. In that case
$\dim_{\R}\Aa= 20$, and $I(h|_{S^2(1)})=1$.\end{exmp}
\begin{exmp}
$h(x_1,x_2,x_3,x_4,x_5)=$$$=(x_1x_2+x_2,3x_3x_5+2x_1,x_1^2+x_2,x_3^2+3x_3,3x_1x_5+x_1,4x_2x_5+x_1,2x_4^2+x_4,x_3^2+x_5)$$
is an immersion on the $4-$dimensional spheres of radius $r=1,10$.
In both cases $\dim_{\R}\Aa= 10$, and $I(h|_{S^4(1)})=0$,
$I(h|_{S^4(10)})=-1$.\end{exmp}
\begin{exmp}
$h(x_1,x_2,x_3,x_4)=$$$=(
x_2x_4+x_4,2x_1x_4+x_3,3x_2x_4+4x_1,3x_3x_4+x_3,x_1x_2+x_3,2x_2x_3+x_1)$$
is an immersion on the $3-$dimensional sphere of radius $r=1$.
Take $u= 3(x_1-y_1)+5(x_2-y_2)-2(x_4-y_4)$. In that case
$\dim_{\R}\Aa=18$, and $I(h|_{S^3(1)})\equiv 1\ \mod 2$.
\end{exmp}


\end{document}